\documentclass[11pt]{article}

\usepackage{euscript}

\usepackage{amsfonts}
\usepackage{amsmath}
\usepackage{amssymb}
\usepackage{latexsym,color}
\usepackage{amsbsy}
\usepackage{graphicx}

\newtheorem{thm}{Theorem}[section]
\newtheorem{con}{Conjecture}[section]
\newtheorem{lemma}[thm]{Lemma}
\newtheorem{cor}[thm]{Corollary}
\newtheorem{pro}[thm]{Proposition}
\newtheorem{example}[thm]{Example}
\newtheorem{definition}[thm]{Definition}

\newtheorem{remark}[thm]{Remark}
\newtheorem{Algorithm}[thm]{Algorithm}

\newcommand{\comment}[1]{}   
 
\newcommand{\CS}{{\cal{X}}}

\newcommand{\Bx}{B_X}

\newcommand{\Co}{{\mathbb C}}   
\newcommand{\N}{{\mathbb N}}   
\newcommand{\Zi}{{\mathbb Z}}   

\newcommand{\pr}{\parallel}

\newcommand{\D}{{\cal{S}}}
\newcommand{\K}{{\cal{K}}}
\newcommand{\I}{{\cal{I}}}
\newcommand{\J}{{\cal{J}}}
\newcommand{\R}{{\cal{R}}}
\newcommand{\M}{{\cal{M}}}
\newcommand{\Z}{{\cal{Z}}}
\newcommand{\C}{{\cal{C}}}

\newcommand{\ul}[1]{{#1}^-}
\newcommand{\ob}[1]{{{\overline{{#1}}}}}

\newcommand{\ncom}{\newcommand}
\ncom{\ns}{\normalsize}
\ncom{\la}{\lambda}
\ncom{\bm}{\boldmath}
\ncom{\noi}{\noindent}
\ncom{\bq}{\begin{equation}}
\ncom{\eq}{\end{equation}}  
\ncom{\beqn}{\begin{eqnarray*}}
\ncom{\eeqn}{\end{eqnarray*}}  
\ncom{\ba}{\begin{array}}

\ncom{\ea}{\end{array}}
\ncom{\beq}{\begin{eqnarray}}
\ncom{\eeq}{\end{eqnarray}}
\ncom{\nno}{\nonumber}
\ncom{\hs}{\mbox{\hspace{.25cm}}}
\ncom{\rar}{\rightarrow}
\ncom{\Rar}{\Rightarrow}  
\ncom{\noin}{\noindent}   
\ncom{\bc}{\begin{center}}
\ncom{\ec}{\end{center}}  
\ncom{\sz}{\scriptsize}   
\ncom{\fpd}{\Phi(\pi^{'})}
\ncom{\fp}{\Phi(\pi) }
\ncom{\nk}{\left< \begin{array}{c}
                       n\\k \end{array} \right>}
\ncom{\nd}{1^{'},2^{'},\cdots,n^{'}}
\ncom{\de}{\bigtriangleup (F_{2n},\leq)}
\ncom{\del}{\bigtriangleup} 
\ncom{\cov}{<\!\!\!\!\cdot }

\ncom{\bt}{\begin{thm}}
\ncom{\bcon}{\begin{con}}
\ncom{\et}{\end{thm}}
\ncom{\econ}{\end{con}}
\ncom{\bl}{\begin{lemma}}
\ncom{\el}{\end{lemma}}  
\ncom{\bco}{\begin{cor}} 
\ncom{\ds}{\displaystyle}
\ncom{\eco}{\end{cor}}   
\ncom{\bp}{\begin{pro}}  
\ncom{\ep}{\end{pro}}    
\ncom{\bex}{\begin{example}}
\ncom{\eex}{\end{example}}  
\ncom{\bd}{\begin{definition}}
\ncom{\ed}{\end{definition}}  
\ncom{\brm}{\begin{remark}}   
\ncom{\erm}{\end{remark}}     
\ncom{\bal}{\begin{Algorithm}}
\ncom{\eal}{\end{Algorithm}}  
\ncom{\ol}{\overline}

\ncom{\pf}{\noi {\bf Proof  }}
\ncom{\be}{\begin{enumerate}} 
\ncom{\ee}{\end{enumerate}}   
\ncom{\s}{\subset}
\ncom{\T}{{\cal T}}
\ncom{\B}{{\cal B}}
\ncom{\A}{{\cal A}}

\title{\Large{{\textcolor{black} {\bf Wreath product action on 
generalized Boolean algebras}}}}

\author{{\textcolor{black} { Ashish Mishra and Murali K. Srinivasan}} \\
{\em  {Department of Mathematics}}\\
{\em  {Indian Institute of Technology, Bombay}}\\
{\em  {Powai, Mumbai 400076, INDIA}}\\
{\bf  \texttt{ashishm@math.iitb.ac.in}}\\ 
{\bf  \texttt{murali.k.srinivasan@gmail.com}}}

\begin{document}
\date{}
\maketitle

\begin{abstract}
Let $G$ be a finite group acting on the finite set $X$ such that the 
corresponding (complex) permutation representation is multiplicity free.
There is a natural rank and order preserving  
action of the wreath product $G\sim S_n$ on the generalized Boolean algebra 
$\Bx(n)$.
We explicitly block diagonalize the commutant 
of this action.
\end{abstract}

\section{Introduction} 

This paper is inspired by the following two results:

(i) Explicit diagonalization of the ``Bose-Mesner algebra" (= commutant of a
certain wreath product action) of the generalized Johnson scheme, by
Ceccherini-Silberstein, Scarabotti, and Tolli {\bf \cite{cst1}}.

(ii) Explicit block diagonalization of the commutant of the wreath product
action on the nonbinary analog of the Boolean algebra, due to Dunkl {\bf
\cite{du1}} and
Gijswijt, Schrijver, and Tanaka {\bf \cite{gst}}.

A natural question suggested by these results is to explicitly block
diagonalize the commutant of the wreath product action on the generalized
Boolean algebra. To do that is one of the aims of this paper. 
Our second aim is to recast the results of {\bf \cite{cst1}} on the generalized
Johnson scheme, presented there
in the language of harmonic analysis (Gelfand pairs, induced
representations, spherical functions etc.), in purely combinatorial terms.
This is achieved by means of the concepts of semisymmetric Jordan basis and
upper Boolean decomposition. These notions are implicit in {\bf
\cite{cst1, gst}} 
and were stated explicitly in {\bf \cite{sr2}}. 
They allow a simple reduction
to the case of symmetric
group action on Boolean algebras, given in {\bf \cite{du1, gst}}.
In the
rest of this introduction we define our objects of study and state our
result.

Let $\mbox{Mat}(n\times n)$ denote the algebra of complex $n\times n$
matrices and let $\A \subseteq \mbox{Mat}(n\times n)$ denote a
$*$-subalgebra. Then the noncommutative analog of the spectral theorem 
asserts that there exists a
{\em{block diagonalization}} of $\A$, i.e.,  there exists a   
$\{1,\ldots ,n\} \times S$ unitary matrix $N$, for some index set $S$ of cardinality
$n$,
and positive integers
$p_0,q_0,\ldots ,p_m,q_m$ such that ${N}^*\A N$ is equal to the
set of all $S\times S$
block-diagonal matrices
\beq \label{bd1} &
\left( \ba{cccc}  C_0 & 0 & \ldots &0\\
            0   & C_1 &\ldots & 0\\
            \vdots & \vdots &\ddots &\vdots\\
            0 & 0 & \ldots & C_m \ea \right), 
&
\eeq
where each $C_k$ is a block-diagonal matrix with $q_k$ repeated, identical
blocks of order $p_k$
\beq \label{bd2} &   
C_k = \left( \ba{cccc}  B_k & 0 & \ldots &0\\
            0   & B_k &\ldots & 0\\
            \vdots & \vdots &\ddots &\vdots\\
            0 & 0 & \ldots & B_k \ea \right). 
&
\eeq
Thus $p_0^2 + \cdots + p_m^2 = \mbox{dim }\A$ and $p_0q_0 + \cdots +
p_mq_m = n$. The numbers $p_0,q_0,\ldots ,p_m,q_m$ and $m$ are uniquely
determined (upto permutation of the indices) by $\A$.

By dropping duplicate blocks we get a $*$-isomorphism 
$$\Phi: \A \cong \bigoplus_{k=0}^m \mbox{Mat}(p_k\times p_k).$$

In an
{\em explicit block diagonalization} we need to know this isomorphism
explicitly, i.e., we need to know  the entries in the
image $\Phi(M)$, for $M$ varying over a suitable basis of $\A$. When $\A$ is
commutative we have $p_k=1$ for all $k$ and we speak of {\em explicit
diagonalization}.

We now define the $*$-algebras to be considered in this paper.
Let $B(n)$ denote the set of all subsets of 
$[n]= \{1,2,\ldots ,n\}$. 
For a finite set $S$, let $V(S)$ denote the complex vector space with $S$ as
basis.

Let $G$ be a finite group acting on the finite set $X$. Assume that the
corresponding permutation representation on $V(X)$ is multiplicity free.
This implies, in
particular, that the action is transitive.

Let $L_0$ be a symbol not in $X$ and let $Y$ denote the {\em alphabet}
$Y=\{L_0\}\cup X$. We call the elements of $X$ the {\em nonzero} letters  in
$Y$.
Define
$ \Bx(n) = \{(a_1,\ldots ,a_n) : a_i \in Y \mbox{ for all }i\}$, the set
of all $n$-tuples of elements of $Y$ (we use $L_0$ instead of $0$ 
for the zero letter for later convenience. We do
not want to confuse the letter $0$ with the vector $0$). 
Given $a = (a_1,\ldots ,a_n)\in \Bx(n)$, define the {\em support} of
$a$ by $S(a) = \{i\in [n] : a_i \not= L_0 \}$. 
For $0\leq i \leq
n$, $\Bx(n)_i$ denotes the set of all elements $a\in \Bx(n)$ with $|S(a)|=i$. 
We have 
$$|\Bx(n)|=(|X|+1)^n,\;\;|\Bx(n)_i|=\binom{n}{i} |X|^i.$$
(We take the binomial coefficient $\binom{n}{k}$ to be 0 if $n<0$ or $k<0$).

Let $S_n$ denote the symmetric group on $n$ letters.
There is a natural action of the wreath product $G\sim S_n$ (see
{\bf\cite{cst2,m}}) on $\Bx(n)$ and
$\Bx(n)_i$:
permute the $n$ coordinates followed by
independently acting on the nonzero letters by elements of $G$. In detail,
given $(g_1,g_2,\ldots ,g_n,\pi)\in G\sim S_n$ (where $\pi\in S_n$ and
$g_i\in G$ for all $i$) and $a=(a_1,\ldots
,a_n)\in \Bx(n)$, we have
$(g_1,\ldots ,g_n,\pi) (a_1,\ldots ,a_n) =
(b_1,\ldots ,b_n)$, where $b_i = g_i a_{\pi^{-1}(i)}$, if $a_{\pi^{-1}(i)}$ is a
nonzero letter and $b_i = L_0$, if $a_{\pi^{-1}(i)}=L_0$.

We give $V(Y)$ and $V(\Bx(n))$ the standard inner products, i.e., we declare
$Y$ and $\Bx(n)$ to be orthonormal bases.
 
We represent elements of $\mbox{End}(V(\Bx(n)))$ (in the standard basis
$\Bx(n)$)
as $\Bx(n)\times \Bx(n)$ matrices (we think of $V(\Bx(n))$ as column vectors
with coordinates indexed by $\Bx(n)$). 
For $a,b\in \Bx(n)$, the entry in row $a$, column $b$ of a matrix $M$ will
be denoted $M(a,b)$. The matrix corresponding to $f\in
\mbox{End}(V(\Bx(n)))$ will be denoted $M_f$. We use similar notations for
$\Bx(n)_i \times \Bx(n)_i$ matrices corresponding to elements of
$\mbox{End}(V(\Bx(n)_i))$.

Set  
\beqn \A_X(n) &=& \{ M_f : f \in \mbox{End}_{G\sim S_n}(V(\Bx(n)))\},\\
 \B_X(n,i) &=& \{ M_f : f \in \mbox{End}_{G\sim S_n} (V(\Bx(n)_i))\}.\eeqn
Thus $\A_X(n)$ and $\B_X(n,i)$ are $*$-algebras of matrices.

In the paper {\bf \cite{cst1}} it is 
shown that the action of $G\sim S_n$ on $V(\Bx(n)_i)$ is
multiplicity free and the commutant $\B_X(n,i)$ of this  action is explicitly
diagonalized. This generalizes the case of the Johnson scheme {\bf\cite{de1}}
(corresponding to $|X|=1$) and the nonbinary Johnson scheme {\bf\cite{tag}}
(corresponding to $G =$ group of all permutations of $X$). For other
examples of this framework see {\bf\cite{cst1,cst2}}.
We consider the related
problem of explicitly block diagonalizing the
commutant $\A_X(n)$ of the action of  $G\sim S_n$ on $V(\Bx(n))$. In the process we
also present an alternative combinatorial approach to the generalized
Johnson scheme of {\bf \cite{cst1}}.

Let $f:V(\Bx(n))\rar V(\Bx(n))$ be linear. Then $f$
is $G\sim S_n$-linear if and only if
\beq \label{gl}
M_f(a,b)&=&M_f(\tau (a),\tau (b)),
\mbox{ for all } a,b\in \Bx(n),\;\tau \in G\sim S_n, \eeq
i.e., $M_f$ is constant on the orbits of the action of $G\sim S_n$ 
on $\Bx(n) \times
\Bx(n)$. We now determine these orbits.

Denote the orbits of $G$ on $X\times X$ by
$$Z_0,Z_1,\ldots ,Z_m,$$
where $Z_0=\{(x,x):x\in X\}$. Note that $m=0$ iff $|X|=1$. Unless otherwise
stated we shall assume $m\geq 1$ (or, equivalently, $|X|\geq 2$). The theory
to be presented in this paper solves the $|X|\geq 2$ case in terms of the
$|X|=1$ case (which is summarized in Section 2 below).

Let
$$C(n,m)=\{(l_1,\ldots ,l_m)\in \N^{m} : l_1 +\cdots +l_m
=n\},$$
where $\N =\{0,1,2,\ldots\}$, denote the set of all {\em weak compositions} of $n$ with $m$ parts. We have
$|C(n,m)|=\binom{n+m-1}{n}=\binom{n+m-1}{m-1}$. 

Given $(a,b)\in \Bx(n)\times \Bx(n)$, where $a=(a_1,\ldots ,a_n)$
and $b=(b_1,\ldots ,b_n)$, define  
$$\C(a,b)=(l_0,\ldots ,l_m)\in C(|S(a)\cap S(b)|,m+1)$$
by 
$$l_j=|\{i\in S(a)\cap S(b) : (a_i,b_i)\in Z_j\}|, \;j=0,\ldots ,m.$$

It is now easily seen that $(a,b),(c,d)\in
\Bx(n) \times \Bx(n)$ are in the same $G\sim S_n$-orbit if and only if
\beq \label{oc}
&|S(a)|=|S(c)|,\;|S(b)|=|S(d)|,\;|S(a)\cap
S(b)|=|S(c)\cap S(d)|,\; \C(a,b) = \C(c,d).
\eeq
For $0\leq i,j,t \leq n$ and $l=(l_0,\ldots ,l_m) \in \N^{m+1}$ 
let $M_{i,j}^{t,l}$ be the $\Bx(n) \times
\Bx(n)$ matrix given by
$$M^{t,l}_{i,j}(a,b) = \left\{ \ba{ll}
                            1 & \mbox{if }|S(a)|=i,\;|S(b)|=j,\;
                                |S(a)\cap S(b)|=t,\; 
                                \C(a,b)=l, \\
                            0 & \mbox{otherwise.}
                            \ea
                    \right.
$$
Define
\beqn \I_X(n) &=& \{(i,j,t,l) : i,j,t \in \N,\; t\leq i,\;t\leq j,\; i+j-t\leq
n,\;l \in C(t,m+1)\}. \eeqn

It follows from (\ref{gl}) and (\ref{oc}) that $\{M_{i,j}^{t,l} :
(i,j,t,l)\in \I_X(n)\}$ is a basis of $\A_X(n)$.
We have 
\beq \label{card1}
&|\I_X(n)|= \mbox{ dim }\A_X(n)
= {{n+m+3}\choose m+3},&\eeq
since $(i,j,t,l)\in \I_X(n)$ if and only if 
$(i-t)+(j-t)+ l_0 + \cdots +l_m\leq n$ and
all $m+3$ terms are nonnegative.

For $0\leq i \leq n$ define
\beqn \I_X(n,i) &=& \{(t,l) : t\in \N,\;t\leq i,\;2i-t\leq n,\;l\in
C(t,m+1)\}.\eeqn
We have that $\{M_{i,i}^{t,l} : (t,l)\in \I_X(n,i)\}$ is a basis of
$\B_X(n,i)$ (here we think of $M_{i,i}^{t,l}$ as $\Bx(n)_i\times \Bx(n)_i$
matrices). We have
\beq \label{card2}
\mbox{ dim }\B_X(n,i) &=& \left\{ \ba{ll}  
\sum_{t=0}^i |C(t,m+1)| = \binom{m+i+1}{m+1}, & \mbox{ if } i\leq n/2,\\&\\
\sum_{t=2i-n}^i |C(t,m+1)| = \binom{m+i+1}{m+1} - \binom{m+2i-n}{m+1}, &
\mbox{ if } i > n/2, \ea \right.
\eeq
where we have used the identity $\sum_{t=0}^i \binom{m+t}{t}
= \binom{m+i+1}{m+1}$ (Exercise 3(a) in Chapter 1 of {\bf\cite{st}}).

Note that if $m=0$ or $1$ then, for $x,y\in X$, $(x,y)$ and $(y,x)$ are in
the same $G$-orbit of $X\times X$. Thus $M_{i,i}^{t,l}$ is symmetric and it
follows that $\B_X(n,i)$ is commutative, since a $*$-subalgebra of
$\mbox{Mat}(n\times n)$ having a basis consisting of symmetric matrices 
is commutative. Hence $V(\Bx(n)_i)$ is a multiplicity free $G\sim
S_n$-module. This argument does not apply when $m\geq 2$.

We now describe our explicit block diagonalization of $\A_X(n)$. We need to
make a number of definitions.

Our poset terminology follows {\bf \cite{st}}.
Let $P$ be a finite {\em graded poset} with 
{\em rank function}
$r: P\rar \N$. The {\em rank} of $P$ is 
$r(P)=\mbox{max}\{r(x): x\in P\}$ and,
for $i=0,1,\ldots ,r(P)$, $P_i$ denotes the set of elements of $P$ of rank
$i$. For a subset $S\subseteq P$, we set $\mbox{rankset}(S) = \{r(x):x\in 
S\}$.

Let $P$ be a graded poset with $n=r(P)$. Then we have
$ V(P)=V(P_0)\oplus V(P_1) \oplus \cdots \oplus V(P_n)$ (vector space direct
sum).
An element $v\in V(P)$ is {\em homogeneous} if $v\in V(P_i)$ for some $i$,
and if $v\not= 0$, we extend the notion of rank to nonzero homogeneous elements by writing
$r(v)=i$. Given an element $v\in V(P)$, write $v=v_0 + \cdots +v_n,\;v_i \in
V(P_i),\;0\leq i \leq n$. We refer to the $v_i$ as the {\em homogeneous
components} of $v$. A subspace $W\subseteq V(P)$ is {\em homogeneous} if it
contains the homogeneous components of each of its elements. For a
homogeneous subspace $W\subseteq V(P)$ we set $\mbox{rankset}(W)=\{r(v) : v
\mbox{ is a nonzero homogeneous element of } W\}$.

The {\em
up operator}  $U:V(P)\rar V(P)$ is defined, for $x\in P$, by
$U(x)= \sum_{y} y$,
where the sum is over all $y$ covering $x$. Similarly,
the {\em
down operator}  $D:V(P)\rar V(P)$ is defined, for $x\in P$, by
$D(x)= \sum_{y} y$,
where the sum is over all $y$ covered by  $x$.

Let $\langle , \rangle$ denote the standard inner product on $V(P)$,
i.e.,
$\langle x,y \rangle =
\delta (x,y)$ (Kronecker delta), for $x,y\in P$.
The {\em length} $\sqrt{\langle v, v \rangle }$ of $v\in V(P)$ is denoted
$\pr v \pr$.

Let $P$ be a graded poset of rank $n$.
A {\em graded Jordan chain} in $V(P)$ is a sequence
\beq \label{gjc}
&c=(v_1,\ldots ,v_h)&
\eeq 
of nonzero homogeneous elements of $V(P)$
such that $U(v_{i-1})=v_i$, for
$i=2,\ldots h$, and $U(v_h)=0$ (note that the
elements of this sequence are linearly independent, being nonzero and of
different ranks). We say that $c$ {\em
starts} at rank $r(v_1)$ and {\em ends} at rank $r(v_h)$.
A {\em graded Jordan basis} of $V(P)$ is a basis of $V(P)$  
consisting of a disjoint union of graded Jordan chains   
in $V(P)$. 

The graded Jordan
chain (\ref{gjc}) is said to be a {\em symmetric Jordan chain} (SJC) if
the sum of the starting and ending ranks of $c$ equals $n$, i.e.,   
$r(v_1) + r(v_h) = n$ 
if $h\geq
2$, or $2r(v_1)= n$ if $h=1$.
A {\em symmetric Jordan basis} (SJB) of $V(P)$ is a basis of $V(P)$
consisting of a disjoint union of symmetric Jordan chains 
in $V(P)$. 

The graded Jordan
chain (\ref{gjc}) is said to be a {\em semisymmetric Jordan chain} (SSJC) if
the sum of the starting and ending ranks of $c$ is $\geq n$.
A {\em semisymmetric Jordan basis} (SSJB) of $V(P)$ is a basis of $V(P)$
consisting of a disjoint union of semisymmetric Jordan chains 
in $V(P)$. 

For $0\leq k \leq n$ we set $\ul{k}=\mbox{ max}(0,2k-n)$ (note that $\ul{k}$
depends on both $k$ and $n$. The $n$ will always be clear from the context).
Note that
$0\leq \ul{k} \leq k$ and $k\leq n + \ul{k} - k$.
For a SSJC $c$ in $V(P)$, starting at rank $i$ and ending at rank $j$,
we define the {\em offset} of $c$ to be $i+j-n$. It is easy to see that if 
an SSJC of offset $s$ starts at rank $k$ then the chain ends at 
rank $n+s-k$ and we have $\ul{k}\leq s \leq k$. Also note that
the conditions $0\leq k \leq n,\;\ul{k}\leq s \leq k$ are easily seen to be
equivalent to 
$0\leq s\leq k \leq \lfloor \frac{n+s}{2} \rfloor \leq n$. 

We use the following terminology in connection with the notion of 
block diagonal form:
let $A$ be a finite set and let $N$ be a $A\times A$ matrix. Given a
partition $A=A_1\cup \cdots \cup A_k$ of $A$ into pairwise disjoint nonempty
subsets $A_1,\ldots ,A_k$ we say that 
{\em $N$ is in block diagonal form with a
block corresponding to each $A_i$} if $A(x,y)=0$ whenever $(x,y)\not\in
(A_1\times A_1) \cup \cdots \cup (A_k \times A_k)$.

Suppose we have an orthogonal graded Jordan basis $O$ of $V(P)$.
Normalize the vectors in $O$ to get an
orthonormal basis $O'$.
Let $(v_1,\ldots ,v_h)$ be a graded Jordan chain in $O$.
Put
$v_u' = \frac{v_u}{\pr v_u \pr}$ and $\alpha_u = \frac{\pr v_{u+1} \pr}{\pr
v_{u} \pr},\;1\leq u \leq h$ (we take $v_0'=v_{h+1}'=0$). 
We have,
for $1\leq u \leq h$,
\beq \label{trick}   
U(v_{u}')=\frac{U(v_{u})}{\pr v_{u} \pr}=\frac{v_{u+1}}{\pr v_{u}
\pr}=\alpha_{u} v_{u+1}'.&&
\eeq
Thus the matrix of $U$ with respect to $O'$ is in block diagonal
form, with a block
corresponding to each (normalized) graded Jordan chain in $O$, and with
the block corresponding to $(v_1',\ldots ,v_h')$ above being a lower
triangular
matrix with subdiagonal $(\alpha_1 ,\ldots ,\alpha_{h-1})$ and $0$'s
elsewhere.

Now note that the matrices, in the standard basis $P$, 
of $U$ and $D$ are real and transposes of
each other. 
Since $O'$ is orthonormal 
with respect to  the standard inner product, it follows that the matrices of
$U$ and $D$, in the basis $O'$, must be adjoints of each other.
Thus the matrix of $D$ with respect to $O'$ is in block diagonal
form, with a block
corresponding to each (normalized) graded Jordan chain in $O$, and with
the block corresponding to $(v_1',\ldots ,v_h')$ above being an upper
triangular
matrix with superdiagonal $(\alpha_1 ,\ldots ,\alpha_{h-1})$ and $0$'s
elsewhere.
Thus, for $0\leq u \leq h-1$,  we have 
\beq \label{trick1} 
D(v_{u+1}')=\alpha_{u} v_{u}'.&&
\eeq
It follows that the subspace spanned by each graded Jordan chain in $O$ is
closed under $U$ and $D$. 

The {\em Boolean algebra} is the rank-$n$ graded poset obtained by partially
ordering $B(n)$ by inclusion. The rank of a subset is given by cardinality.

Given $a = (a_1,\ldots ,a_n), b=(b_1,\ldots ,b_n)\in \Bx(n)$, 
define $a \leq b$ 
provided $S(a) \subseteq S(b)$ and $a_i = b_i$ for all
$i\in S(a)$. It is easy to see that this makes $\Bx(n)$ into a
rank-$n$ graded
poset with rank of $a$ given by $|S(a)|$. 
We call $\Bx(n)$ a {\em generalized
Boolean algebra}. 
Clearly, when $|X|=1$,
$\Bx(n)$ is order isomorphic to $B(n)$. When $G$ is the group of all
permutations of $X$ 
we have two orbits for the action on $G$ on
$X\times X$ (i.e., $m=1$) and we call $\Bx(n)$ the {\em nonbinary analog} of
the Boolean algebra.
We use $U_n$
to denote the up operator on both the posets 
$V(B(n))$ and $V(\Bx(n))$. Note that the action of $G\sim S_n$ on $\Bx(n)$
is order and rank preserving and that $U_n$ is $G\sim S_n$-linear.
Also note that the inner product on $V(\Bx(n))$ is $G\sim S_n$-invariant.

Consider the permutation representaion of $G$ on $V(X)$. A map $f:V(X)\rar
V(X)$ is $G$-linear iff the $X\times X$ matrix $M_f$ representing $f$ (in
the standard basis $X$) is constant on the orbits of the action of $G$ on
$X\times X$. It follows that $\mbox{dim }\mbox{End}_G(V(X))=m+1$.

Consider the canonical orthogonal decomposition
$$V(X)=W_0\oplus\cdots \oplus W_m,$$
of $V(X)$ into distinct irreducible $G$-submodules, where $W_0$ corresponds to the
trivial representation. Set $d_i=\mbox{dim }W_i$, $i=0,1,\ldots ,m$, so
$d_0=1$ and $d_0+\cdots +d_m = |X|$.

For $u=0,\ldots ,m$, define the linear map $f_u:V(X)\rar V(X)$ by 
$$f_u(y)=\sum x,\;\;\;y\in X,$$
where the sum is over all $x\in X$ with $(x,y)\in Z_u$. The matrix of $f_u$,
in the standard basis $X$, is the $0-1$ matrix with a 1 in position $(x,y)$
iff $(x,y)\in Z_u$. Thus $f_u$ is $G$-linear. Since $W_0,\ldots ,W_m$
are distinct irreducibles we have 
by Schur's lemma that each of $W_0,\ldots ,W_m$ is $f_u$ invariant for each
$u$ and that  the action of each $f_u$ on each of $W_0,\ldots ,W_m$ 
is multiplication by a scalar. For $u,w=0,\ldots , m$, let the action of
$f_u$ on $W_w$ be multiplication by the scalar $\lambda(u,w)$.

Set
\beqn \J_X(n) &=&
\{(k,s,p):0\leq k \leq n,\; \ul{k}\leq s \leq k, \;p\in C(s,m)\}\\
&=& \{(k,s,p):0\leq s \leq n,\; s\leq k \leq \lfloor\frac{n+s}{2}\rfloor, 
\;p\in C(s,m)\}.
\eeqn

For $(k,s,p)\in \J_X(n)$ with $p=(p_1,\ldots ,p_m)$ define
\beqn 
\mu(n,k,s,p) 
&=& \binom{n}{n-s,p_1,\ldots ,p_m}\;
d_1^{p_1}\ldots d_m^{p_m}
\left\{ {{n-s}\choose {k-s}}-{{n-s}\choose{k-s-1}}\right\}.
\eeqn
(We take the multinomial coefficient $\binom{n}{k_1,\ldots ,k_r}$ to be zero
if any of $n,k_1,\ldots ,k_r$ is negative or if $n\not=k_1+\cdots +k_r$).

The following two definitions are from {\bf\cite{s}} and {\bf\cite{cst1}}
respectively.

For $i,j,k,t \in \{0,\ldots ,n\}$ define
\beqn
\beta_{i,j,k}^{n,t} & = &\sum_{u=0}^n (-1)^{u-t} \;
\binom{u}{t}\binom{n-2k}{u-k}\binom{n-k-u}{i-u}\binom{n-k-u}{j-u}.
\eeqn

For $l=(l_0,\ldots ,l_m), p=(p_0,\ldots ,p_m) \in \Zi^{m+1}$ define
$$\Lambda(\lambda,l,p) =  \sum 
\left\{ \prod_{w=0}^m \binom{p_w}
{r(0,w),\ldots ,r(m,w)}\right\}
\left\{ \prod_{u=0}^m \prod_{w=0}^m
\lambda(u,w)^{r(u,w)} \right\}, $$   
where the sum is over all $\{0,1,\ldots ,m\}\times \{0,1,\ldots ,m\}$
nonnegative integer matrices $(r(u,w))$ with row sums $l_0,\ldots ,l_m$ 
and column sums $p_0,\ldots ,p_m$. We take the empty sum to be zero. 

For $i,j,t,k,s \in \{0,\ldots ,n\}$, $l\in C(t,m+1)$, and 
$p=(p_1,\ldots ,p_m)\in C(s,m)$,
set $p_0=t-s$, $p^+ =(p_0,p_1,\ldots ,p_m)$ 
and define
\beqn
\alpha_{i,j,k,s}^{n,t,l,p} &=&  (|X|)^{\frac{i+j}{2} -t}\;
                \Lambda(\lambda,l,p^+)\;
                \beta_{i-s,j-s,k-s}^{n-s,t-s}.
\eeqn

For $0\leq k\leq n,\;\ul{k}\leq s \leq k$, and $k\leq i,j \leq n+s-k$, 
define $E_{i,j,k,s}$ to be the $\{k,k+1,\ldots ,n+s-k\}\times \{k,k+1,\ldots
,n+s-k\}$ matrix
with entry in row $i$ and column $j$ equal to 1 and all other
entries 0.

We now state our result. Part (i) is proved in Lemma \ref{mt2},  
part (ii) is proved in Theorem \ref{mt3},  and
parts (iii) and (iv) are proved in  Section 4.
\bt  \label{ebdgba} 
Let $G$ be a finite group acting on the finite set $X$ such that the
corresponding complex permutation representation is multiplicity free
with $m+1$ distinct irreducibles occuring in it. Assume $m\geq 1$. Then

(i) There exists an orthogonal semisymmetric Jordan basis $J_X(n)$ of
$V(\Bx(n))$.

(ii) $J_X(n)$ can be expressed as a disjoint union
$$ J_X(n)= \cup_{(k,s,p)\in \J_X(n)} J_X(n,k,s,p),$$
where $J_X(n,k,s,p)$ consists of $\mu(n,k,s,p)$ semisymmetric Jordan chains
starting at rank $k$ and ending at rank $n+s-k$.

Let $J_X(n)'$ denote the orthonormal basis obtained by normalizing the
vectors in $J_X(n)$. Define a  
$\Bx(n) \times J_X'(n)$ unitary matrix $M(n)$ as follows: for $v\in J_X'(n)$, 
the column of
$M(n)$ indexed by $v$ is the coordinate vector of $v$ in the standard  
basis $\Bx(n)$. 
 
(iii) $M(n)^* \A_X(n) M(n)$ consists of all $J_X'(n) \times J_X'(n)$  
block diagonal matrices with a block
corresponding to each (normalized) SSJC in $J_X(n)$ and such that, 
for each $(k,s,p)\in \J_X(n)$, the $\mu(n,k,s,p)$ blocks corresponding
to the SSJC's in $J_X(n,k,s,p)$ are identical.

(iv) Conjugating by $M(n)$ and, for each $(k,s,p)\in \J_X(n)$,
preserving only one among the duplicate
blocks corresponding to the SSJC's in $J_X(n,k,s,p)$,
thus gives a
$*$-algebra isomorphism 
$$\Phi: \A_X(n) \cong \bigoplus_{(k,s,p)\in \J_X(n)} 
\mbox{Mat}((n+s-2k+1)\times (n+s-2k+1)).$$
It will be convenient to re-index
the rows and columns of a block corresponding to a SSJC starting at rank $k$
and having offset $s$ by the set $\{k,k+1,\ldots ,n+s-k\}$.
 
Let $(i,j,t,l) \in \I_X(n)$. Write
$$\Phi (M_{i,j}^{t,l}) = (N_{k,s,p}),\;(k,s,p)\in \J_X(n).$$
Then
$$N_{k,s,p} = \left\{ \ba{ll}
                {{n+s-2k}\choose{i-k}}^{-\frac{1}{2}}\;
                {{n+s-2k}\choose{j-k}}^{-\frac{1}{2}}\;
                \alpha_{i,j,k,s}^{n,t,l,p}\;
                E_{i,j,k,s} & \mbox{if } k\leq i,j \leq n+s-k, \\
                  0    & \mbox{otherwise}.
                 \ea
         \right.
$$ 
\et

\section{Symmetric group action on Boolean algebras}

In this section we summarize results on the $S_n$ action on $B(n)$ in the
form needed in this paper. The main sources are 
Dunkl {\bf\cite{du1}} and Schrijver {\bf\cite{s}}. 
For $0\leq i,j,t \leq n$  
let $M_{i,j}^{t}$ be the $B(n) \times B(n)$ matrix given by
$$M^{t}_{i,j}(a,b) = \left\{ \ba{ll}
                            1 & \mbox{if }|a|=i,\;|b|=j,\;
                                |a\cap b|=t, \\
                            0 & \mbox{otherwise.}
                            \ea
                    \right.
$$
It is clear that $\{ M_{i,j}^t : i,j,t\in \N,\; t\leq i, \;t\leq j,\; i+j-t \leq n\}$ is a
basis of 
$$\A(n) = \mbox{ End}_{S_n}(V(B(n))).$$

The following formulation is 
motivated by the one given in {\bf\cite{s}} (also see {\bf\cite{sr1}}). We
do not use part (v) (a classical result of Delsarte {\bf\cite{de1}}) in this
paper but we have included it for completeness.  

\bt \label{schrijver}
(i) There exists an orthogonal SJB $J(n)$ of $V(B(n))$.

For $0\leq k \leq n/2$, let $J(n,k)$ be the set of all SJC's in $J(n)$
starting at rank $k$. For $0\leq i \leq n$ and $0\leq k \leq
\mbox{min}\{i,n-i\}$ set
$$J(n,i,k) = \{v\in J(n,k) : r(v)=i\}$$
and let $V(n,i,k)$ be the subspace of $V(B(n)_i)$ spanned by $J(n,i,k)$.

(ii) Let $0\leq i \leq n$. Then $V(B(n)_i)$ is a multiplicity free
$S_n$-module and
\beqn V(B(n)_i) &=& \bigoplus_{0\leq k \leq \mbox{min}\{i,n-i\}} V(n,i,k)\eeqn
is the orthogonal decomposition of $V(B(n)_i)$ into distinct irreducible
submodules.

Let $J'(n)$ denote the orthonormal basis of $V(B(n))$ obtained by
normalizing  $J(n)$.  
Define a  
$B(n) \times J'(n)$ unitary matrix $N(n)$ as follows: for $v\in J'(n)$, 
the column of
$N(n)$ indexed by $v$ is the coordinate vector of $v$ in the standard  
basis $B(n)$.  
 
(iii) $N(n)^* \A(n) N(n)$ consists of all $J'(n) \times J'(n)$  
block diagonal matrices with a block
corresponding to each (normalized) SJC in $J(n)$ and such that 
any two SJC's
starting and ending at the same rank give rise to identical blocks.
Thus there are $\binom{n}{k}- \binom{n}{k-1}$ identical
blocks of size $(n-2k+1)\times (n-2k+1)$, 
for $k=0,\ldots , \lfloor n/2 \rfloor$. 

(iv) Conjugating by $N(n)$ and dropping duplicate blocks thus gives a 
$*$-algebra isomorphism 
$$\Phi: \A(n) \cong \bigoplus_{k=0}^{\lfloor n/2 \rfloor}
\mbox{Mat}((n-2k+1)\times (n-2k+1)).$$
It will be convenient to re-index
the rows and columns of a block corresponding to a SJC starting at rank $k$
and ending at rank $n-k$ by the set $\{k,k+1,\ldots ,n-k\}$. 

Let $i,j,t \in \N, t\leq i,\;t\leq j,\;i+j-t\leq n$. Write
$$\Phi (M_{i,j}^t) = (N_0,\ldots ,N_{\lfloor n/2 \rfloor}).$$
Then, for $0\leq k \leq \lfloor n/2 \rfloor$,
$$N_k = \left\{ \ba{ll}
                {\binom{n-2k}{i-k}}^{-\frac{1}{2}}
                {\binom{n-2k}{j-k}}^{-\frac{1}{2}}\;
                \beta_{i,j,k}^{n,t}\;
                E_{i,j,k} & \mbox{if } k\leq i,j \leq n-k, \\
                  0    & \mbox{otherwise}.
                 \ea
         \right.
$$

(v) Let $0\leq i \leq n$,
$\ul{i} \leq t \leq i$, and $0\leq k \leq \mbox{min}\{i,n-i\}$. Substituting
$j=i$ in part (iv) and simplifying we get that the eigenvalue of 
$M_{i,i}^t$ on $V(n,i,k)$ is
$$\sum_{u=0}^n 
(-1)^{u-t}\;\binom{u}{t}\binom{n-k-u}{i-u}\binom{i-k}{i-u}.$$
\et

For a subspace analog (or $q$-analog) of Theorem \ref{schrijver}, see {\bf\cite{bvp,sr2}}. For
explicit constructions of orthogonal SJB's of $V(B(n))$ and $V(B_q(n))$ (the
subspace analog of $V(B(n))$) see {\bf\cite{sr1,sr3}}.

\section{Upper Boolean decomposition}

Let $(V,f)$ be a pair consisting of a finite dimensional inner product space
$V$ (over $\Co$) and a linear operator $f$ on $V$. Let $(W,g)$ be
another such pair. By an isomorphism of pairs $(V,f)$ and $(W,g)$ we mean a
linear isometry (i.e, an inner product preserving isomorphism) $\theta : V 
\rar W$ such that $\theta(f(v)) = g(\theta(v)),\;v\in V$.

Consider the inner product space $V(\Bx(n))$.
An {\em upper Boolean subspace} of rank $t$ is a homogeneous
subspace $W\subseteq V(\Bx(n))$ such that
$\mbox{rankset}(W)=\{n-t,n-t+1,\ldots ,n\}$, $W$ is closed under the up
operator $U_n$, and there is an isomorphism of pairs
$(V(B(t)),\sqrt{|X|}\,U_t)
\cong (W,U_n)$ that sends homogeneous elements to homogeneous elements and   
increases rank by $n-t$.

Consider the identity
\beq \nonumber
(|X|+1)^n &=& (d_1+d_2+\cdots +d_m +2)^n \\ \label{multid}
          &=& \sum_{(l_0,\ldots ,l_m)\in C(n,m+1)} 
              {\binom{n}{l_0,\ldots ,l_m}}\;d_1^{l_1}\;d_2^{l_2}\cdots
d_m^{l_m}\;2^{l_0}.
\eeq

We shall now give a linear algebraic interpretation to the identity above.

Consider the inner
product space $V(Y)$, with $Y$ as an orthonormal basis. 
Make the tensor product
$$ \otimes_{i=1}^n V(Y) = V(Y)\otimes \cdots \otimes V(Y)\;\;(n
\mbox{ factors})$$
into an inner product space by defining 
\beq \label{ipl}
\langle v_1 \otimes \cdots \otimes
v_n , u_1 \otimes \cdots \otimes u_n \rangle &=& \langle v_1 , u_1 \rangle
\cdots \langle v_n , u_n \rangle .
\eeq
There is an isometry  
\beq \label{tp}
V(\Bx(n)) &\cong &
\otimes_{i=1}^n V(Y)
\eeq
given by 
$a=(a_1,\ldots ,a_n) \mapsto {\ob{a}} = a_1 \otimes
\cdots \otimes a_n,\;a\in \Bx(n)$.
The rank function (on nonzero homogeneous elements) 
and the up and down operators, $U_n$ and $D_n$, 
on $V(\Bx(n))$ are transferred to
$\otimes_{i=1}^n V(Y)$ via the isomorphism above. From now onwards, we shall
not distinguish between $V(\Bx(n))$ and $\otimes_{i=1}^n V(Y)$.

Choose orthonormal bases $\B_i$ for $W_i$, $i=0,\ldots ,m$ and set $\B=\B_0
\cup \cdots \cup \B_m$. We assume that $\B_0$ consists of the vector
$$\frac{1}{\sqrt{|X|}}\sum_{x\in X}x = z.$$
Note that orthonormality of $\B$
implies that if $v=\sum_{x\in X}\alpha_x x \in \B_i,\;i=1,\ldots ,m$, then
$\sum_{x\in X}\alpha_x =0$.

We have in $V(Y)$ 
\beq \label{ipt1}
&&U_1(v) = D_1(v) = 0,\; v\in \B_i, \;i=1,\ldots ,m,\\
\label{ipt2} && U_1(z)= D_1(L_0) = 0,\\
\label{ipt3} && U_1(L_0) = \sqrt{|X|} \;z, \;D_1(z)= \sqrt{|X|} \;L_0.
\eeq 
Given $l=(l_0,\ldots ,l_m)\in C(n,m+1)$ set
\beqn 
S(n,m,l)&=& \{ (A_1,\ldots ,A_m) : 
A_i \subseteq [n] \mbox{ for all }i, A_i\cap A_j = \emptyset 
\mbox { for }i\not= j, \mbox{ and } |A_i|=l_i \mbox{ for all $i$}\}. \eeqn

Given $\A=(A_1,\ldots ,A_m)\in S(n,m,l)$ we set
\beqn
\Sigma(\A)&=& A_1\cup \cdots \cup A_m, \\
\M(\A)&=& \{ f=(f_1,\ldots ,f_m) \;|\; f_i :A_i \rar \B_i \mbox{ for all $i$}\}.
\eeqn 
Given $l=(l_0,\ldots ,l_m)\in C(n,m+1)$ define
\beq 
\D_X(n,l)&=& \{ (\A,f) : \A\in S(n,m,l),\;f\in \M(\A)\}\\ 
\K_X(n,l)&=& \{ (\A,f,B) : \A\in S(n,m,l),\;f\in \M(\A),\;B\subseteq
[n] - \Sigma(\A)\}. 
\eeq
Note that
\beqn
|\D_X(n,l)| &=&  {\binom{n}{l_0,\ldots ,l_{m}}}\;d_1^{l_1}\;d_2^{l_2}\cdots
d_m^{l_m},\\
|\K_X(n,l)| &=&  {\binom{n}{l_0,\ldots ,l_m}}\;d_1^{l_1}\;d_2^{l_2}\cdots
d_m^{l_m}\;2^{l_0}.
\eeqn

For $(\A,f,B) \in \K_X(n,l)$, with $\A=(A_1,\ldots ,A_m)$, $f=(f_1,\ldots
,f_m)$, define a vector 
$$v(\A,f,B) = v_1 \otimes \cdots \otimes v_n \in
\otimes_{i=1}^n V(Y)$$ 
by
\beqn
v_i &=& \left\{ \ba{ll}
          f_j(i)& \mbox{if $i\in A_j,$}\\
                                 & \\ 
          z & \mbox{if $i\in B,$}\\
                    & \\
      L_0 & \mbox{if $i\in [n] -(\Sigma(\A) \cup B).$}
                   \ea \right.
\eeqn
Note that $v(\A,f,B)$ is a homogeneous vector in $\otimes_{i=1}^n V(Y)$ 
of rank $|\Sigma(\A)| + |B|$.

Given $l\in C(n,m+1)$ and
$(\A,f)\in \D_X(n,l)$, set
$$K_X(l,\A,f) = \{v(\A,f,B) : B \subseteq ([n] - \Sigma(\A)) \},$$
and define $V_{(l,\A,f)}$ to be 
the subspace of $\otimes_{i=1}^n V(Y)$ spanned by
$K_X(l,\A,f)$. Set 
\beqn K_X(n) &=& \cup_l  \cup_{(\A,f)} K_X(l,\A,f),\eeqn
where the union is over $l\in C(n,m+1)$ and $(\A,f) \in \D_X(n,l)$.

We have, using (\ref{ipt1}), (\ref{ipt2}), and (\ref{ipt3}), the following
formula in $\otimes_{i=1}^n V(Y)$: for $v(\A,f,B)\in K_X(n)$
\beq \label{ubs}
U_n(v(\A,f,B)) &=& \sqrt{|X|} \left\{ \sum_{B'} v(\A,f,B') \right\},
\eeq
where the sum is over all $B' \subseteq ([n] - \Sigma(\A))$ covering $B$.

It follows from the orthonormality of $\B$ and the 
inner product formula (\ref{ipl}) that 
\beq \label{tpf}
\langle v(\A,f,B) , v(\A',f',B') \rangle &=& \delta((\A,f,B),(\A',f',B')),
\eeq
where $v(\A,f,B),\;v(\A',f',B') \in K_X(n)$ and $\delta$ is the Kronecker
delta.

We can summarize the discussion above in the following result. 
\bt \label{oubdt}
(i) For $l\in C(n,m+1),\;(\A,f)\in \D_X(n,l)$, $V_{(l,\A,f)}$ 
is an upper Boolean subspace 
of $\otimes_{i=1}^n V(Y)$ of rank $n - |\Sigma(\A)|$ and 
with orthonormal basis $K_X(l,\A,f)$.
 
(ii) $K_X(n)$ is an orthonormal basis of $\otimes_{i=1}^n V(Y)$.

(iii) We have the following orthogonal decomposition into upper Boolean
subspaces:
\beq \label{ubd}
\bigotimes_{i=1}^n V(Y) &=& \bigoplus_l \bigoplus_{\A} 
\bigoplus_f V_{(l,\A,f)},
\eeq
where the sum is over all $l\in C(n,m+1)$, $\A\in S(n,m,l)$ and
$f\in\M(\A)$.

For each $(l_0,\ldots ,l_m)\in C(n,m+1)$ the right hand side of
(\ref{ubd}) contains $\binom{n}{l_0,\ldots ,l_m}\;d_1^{l_1}\;d_2^{l_2}\cdots
d_m^{l_m}$ upper Boolean subspaces of
rank $l_0$.

\et

Taking dimensions on both sides of (\ref{ubd}) yields a linear algebraic
interpretation to identity (\ref{multid}) above.
Certain problems on the generalized Boolean algebra $\Bx(n)$ 
can be reduced to  
corresponding problems on the Boolean algebra $B(n)$ via
the basis $K_X(n)$. 

\bl \label{mt2}
(i) For $l=(l_0,\ldots ,l_m) \in C(n,m+1)$, $(\A,f)\in \D_X(n,l)$
there exists an orthogonal SSJB $J_X(l,\A,f)$ of $V_{(l,\A,f)}$.

Each SSJC in $J_X(l,\A,f)$ has offset $s=l_1+\cdots +l_m$. 
For $s \leq k \leq \lfloor \frac{n+s}{2} \rfloor$, denote by $J_X(k,l,\A,f)$
the set of all SSJC's in $J_X(l,\A,f)$ starting at rank $k$. Then
$$|J_X(k,l,\A,f)| = \binom{n-s}{k-s} -
\binom{n-s}{k-s-1}.$$

(ii) $J_X(n) = \cup J_X(l,\A,f)$, where the union is over all $l\in
C(n,m+1)$ and $(\A,f)\in \D_X(n,l)$, is an orthogonal
SSJB of $\otimes_{i=1}^n V(Y)$.

\el
\pf
(i) Consider the upper Boolean subspace $V_{(l,\A,f)}$ of rank $l_0$.
Let $\gamma : \{1,2,\ldots ,l_0\} \rar [n] - \Sigma(\A)$ be the unique order
preserving bijection, i.e., $\gamma (i) = i^{th}$ smallest element of
$[n]- \Sigma(\A)$. Denote by $\gamma' : V(B(l_0)) \rar V_{(l,\A,f)}$ 
the linear isometry 
given by $\gamma' (X) = v(\A,f,\gamma (X)),\;X\in B(l_0)$.

Use Theorem \ref{schrijver} (i) to get an orthogonal SJB $J(l_0)$ of $V(B(l_0))$
with respect to $\sqrt{|X|}\;U_{l_0}$ 
(rather than just $U_{l_0}$) and transfer it to $V_{(l,\A,f)}$ via $\gamma'$.
Each SJC in $J(l_0)$ will get transferred
to a SSJC in $V_{(l,\A,f)}$ of offset $s$. 
The number of these SSJC's in
$V_{(l,\A,f)}$
starting at rank $k$ is 
${{n-s}\choose {k-s}} - {{n-s}\choose {k-s-1}}$. 

(ii) This is clear. $\Box$

We shall now make three observations on the action of $G\sim S_n$ on the
basis elements in $J_X(k,l,\A,f)$.

Let $l=(l_0,l_1,\ldots ,l_m)\in C(n,m+1)$ and 
$(\A,f)\in \D_X(n,l)$. Write $\A=(A_1,\ldots ,A_m)$ and
$f=(f_1,\ldots ,f_m)$. 

(a) Let $v\in J_X(k,l,\A,f)$ (i.e., $v$ is a vector lying in one of the
SSJC's in $J_X(l,\A,f)$ starting at rank $k$) of rank $i$ and 
let $\tau=(1,\ldots ,1,\pi)\in G\sim S_n$
($1 = $ identity in $G$). 

Write $\pi \A=(\pi A_1,\ldots ,\pi A_m)=\A'$ and 
$f\pi^{-1} = (f_1\pi^{-1}|_{\pi A_1},\ldots ,f_m \pi^{-1}|_{\pi A_m})=f'$.
We want to show that $\tau v$ is in the subspace spanned by all 
elements of $J_X(k,l,\A',f')$ of rank $i$.

Let $\sigma: [n]-\Sigma(\A)\rar [n]-\Sigma(\A')$ be the unique order
preserving bijection.
There is an obvious isomorphism of pairs
$\Gamma : (V_{(l,\A,f)},U_n) \rar (V_{(l,\A',f')},U_n)$
taking $v(\A,f,B)\in V_{(l,\A,f)}$
to $\Gamma(v(\A,f,B))=v(\A',f',\sigma B)$.

Clearly we have
$\Gamma (v) \in J_X(k,l,\A',f')$ with $r(\Gamma(v))=i$. 
Define a permutation $\pi':[n]\rar [n]$ as
follows: $\pi'(i) = i$ for $i \in \Sigma(\A')$ and $\pi'(i)= \pi
\sigma^{-1} (i)$ for $i\in [n]- \Sigma(\A')$. Set
$\tau'=(1,\ldots ,1,\pi')\in G\sim S_n$. A little reflection shows that 
$\tau v = \tau' \Gamma(v)$. 

Now, by Theorem \ref{schrijver} (ii), the
subspace spanned by all elements of $J_X(k,l,\A',f')$ of rank $i$ is closed
under permutations of $[n]$ that fix $\Sigma(\A')$ and thus  
it follows that $\tau v$ is in this subspace.

(b) Let $v\in J_X(k,l,\A,f)$, $i\in A_j$, $g\in G$, 
and $\tau=(1,\ldots ,1, g,1, \ldots ,1,\mbox{id})\in G\sim S_n$, with
$g$ in the $i$th spot and id denoting the identity of $S_n$.
Then we have 
$$\tau v = \sum_e \alpha_e v_e,$$
where the $\alpha_e$ are scalars, $v_e\in J_X(k,l,\A,e)$, and 
where the sum is over all 
$$e=(f_1,\ldots ,f_{j-1},e_j,f_{j+1},\ldots ,f_m)$$
with $e_j : A_j \rar \B_j$ satisfying $e_j (w) = f_j (w)$, for all $w\in A_j -
\{i\}$.

(c) Let $v\in J_X(k,l,\A,f)$, $i\in [n]-\Sigma(\A)$, $g\in G$, 
and $\tau=(1,\ldots ,1, g,1, \ldots ,1,\mbox{id})\in G\sim S_n$, with
$g$ in the $i$th spot and id denoting the identity of $S_n$. Then we have
$\tau v= v$.

\bt \label{mt3} 
Let $0\leq s \leq n$, $p=(p_1,\ldots ,p_m)\in C(s,m)$ and set
$l=(n-s,p_1,\ldots ,p_m)\in C(n,m+1)$. 

Define
$$V_{(s,p)} = \bigoplus_{\A} \;\bigoplus_{f} V_{(l,\A,f)},\;\;\;
J_X(n,s,p) = \ds{\cup_{(\A,f)}}\; J_X(l,\A,f),$$
where the sum and the disjoint union are over all $(\A,f)\in \D_X(n,l)$.

(i) We have the following orthogonal decomposition
\beq \label{ubd1}
\bigotimes_{i=1}^n V(Y) &=& \bigoplus_{0\leq s \leq n}
\;\;\bigoplus_{p\in C(s,m)} V_{(s,p)}.
\eeq

(ii) $J_X(n,s,p)$ is an orthogonal SSJB of $V_{(s,p)}$ with each SSJC in 
$J_X(n,s,p)$ having offset $s$. For 
$s\leq k \leq \lfloor \frac{n+s}{2} \rfloor $ define
$J_X(n,k,s,p)$ to be the collection of all chains in $J_X(n,s,p)$ starting
at rank $k$.
Then $J_X(n,k,s,p)$ consists of $\mu(n,k,s,p)$ semisymmetric Jordan 
chains starting
at rank $k$ and ending at rank $n+s-k$. We have
$$J_X(n) = \cup_{(k,s,p)\in \J_X(n)} J_X(n,k,s,p).$$
\et

\pf Clear. $\Box$

For $0\leq i \leq n$, define
\beqn \J_X(n,i) &=& 
\{ (k,s,p) \in \J_X(n) : k\leq i \leq n+s-k\}\\
&=& \{ (k,s,p) \in \J_X(n) : s\leq k \leq \mbox{ min}\{i,n+s-i\}\}.
\eeqn
For $0\leq i \leq n$ and $(k,s,p)\in \J_X(n,i)$, define
$$J_X(n,i,k,s,p) = \{ v\in J_X(n,k,s,p) : r(v) = i\}$$
and $V_X(n,i,k,s,p)$ to be the subspace of $V(\Bx(n)_i)$ spanned by
$J_X(n,i,k,s,p)$. We have the following orthogonal decomposition
\beq \label{odgjs}
V(\Bx(n)_i) &=& \bigoplus_{(k,s,p)\in \J_X(n,i)} V_X(n,i,k,s,p).
\eeq 

\bt \label{edgjs} 
(Ceccherini-Silberstein,
Scarabotti, and Tolli {\bf\cite{cst1}})
The action of $G\sim S_n$ on $V(\Bx(n)_i)$ is
multiplicity free  and (\ref{odgjs}) gives the decomposition into
distinct irreducible submodules.
\et
\pf Let $0\leq i \leq n$ and
$(k,s,p)\in \J_X(n,i)$. It follows from the three properties (a), (b), and (c)
stated above
that $V_X(n,i,k,s, p)$ is a $G\sim S_n$-submodule of $\otimes_{i=1}^n V(Y)$.
   
The decomposition (\ref{odgjs}) is indexed by $\J_X(n,i)$ and a basis of
$\mbox{End}_{G\sim S_n}(V(\Bx(n)_i))$ is indexed by $\I_X(n,i)$. The result will follow if we show
that $\J_X(n,i)$ and $\I_X(n,i)$ have the same cardinality.

\noi Case (i) $i\leq n/2$: We have
\beqn
\I_X(n,i) &=& \{ (t,l) : 0\leq t \leq i,\; l\in C(t,m+1)\}\\     
\J_X(n,i) &=& \{ (k,s,p) : 0\leq s\leq k \leq i,\; p\in C(s,m)\}\\     
          &=& \{ (k,p) : 0 \leq k \leq i,\; p\in C(k,m+1)\}
\eeqn
Clearly $\I_X(n,i)$ and $\J_X(n,i)$ have the same cardinality.     

\noi Case (ii) $i > n/2$: We have
\beqn
\I_X(n,i) &=& \{ (t,l) : 2i-n\leq t \leq i,\; l\in C(t,m+1)\}\\
          &=& \{ (t,s,p): 2i-n \leq t \leq i,\;0\leq s \leq t,\;p\in
C(s,m)\}\\    
          &=& \{ (t,s,p): 0\leq s ,\;\mbox{ max}\{2i-n,s\}\leq t\leq
i,\;p\in C(s,m)\}.\\    
\J_X(n,i) &=& \{ (k,s,p) : 0\leq s,\;s\leq k \leq \mbox{
min}\{i,n+s-i\},\; p\in C(s,m)\}.     
\eeqn
Clearly $\I_X(n,i)$ and $\J_X(n,i)$ have the same cardinality. $\Box$     

We have from Theorem \ref{edgjs} above that, for $0\leq i \leq n$, 
the $V_X(n,i,k,s,p),\;(k,s,p)\in \J_X(n,i)$
are the common eigenspaces for $M_{i,i}^{t,l},\;(t,l)\in \I_X(n,i)$. The
eigenvalues will follow from part (iv) of Theorem \ref{ebdgba} (see Theorem
\ref{ev}).

\section{Explicit block diagonalization }

We shall need the following result.
\bl  \label{bi} For $m\geq 1$ we have
\beqn 
\sum_{k=0}^n \sum_{s=\ul{k}}^k\;(n+s-2k+1)^2\;\binom{s+m-1}{m-1} &=&
\binom{n+m+3}{m+3}. \eeqn
\el
\pf We shall use the following well known identities:

(i) For $n$ odd, $1^2 + 3^2 + 5^2 + \cdots + n^2 = \binom{n+2}{3}$.

(ii) For $n$ even, $2^2 + 4^2 + 6^2 + \cdots + n^2 = \binom{n+2}{3}$.

(iii) For $r,s\geq 0$ we have
$$ \sum_{k=0}^n \binom{k+r}{r}\;\binom{n-k+s}{s} = \binom{n+r+s+1}{r+s+1}.$$
(For the proof multiply the expansions 
$\sum_{k\geq 0} \binom{k+r}{r}x^k=\frac{1}{(1-x)^{r+1}}$ and
$\sum_{k\geq 0} \binom{k+s}{s}x^k=\frac{1}{(1-x)^{s+1}}$ and compare
coefficients of $x^n$ on both sides.)

We have
\beqn
\lefteqn{\sum_{s=0}^n \sum_{k=s}^{\lfloor \frac{n+s}{2} \rfloor}
(n+s-2k+1)^2 \binom{s+m-1}{m-1}}\\ 
&=& 
\sum_{s=0,\; n+s \mbox{ \scriptsize{even}}}^n (1^2+3^2+\cdots +
(n-s+1)^2)\binom{s+m-1}{m-1} \\ && +
\sum_{s=0,\; n+s \mbox{ \scriptsize{odd}}}^n (2^2+4^2+\cdots +
(n-s+1)^2)\binom{s+m-1}{m-1}\\
&=& \sum_{s=0}^n \binom{n-s+3}{3}\binom{s+m-1}{m-1}\\
&=& \binom{n+m+3}{m+3}.\; \Box
\eeqn

Consider the linear operator on $V(\Bx(n))$
whose matrix with respect to the standard basis $\Bx(n)$ 
is $M_{i,j}^{t,l}$. Transfer this operator to 
$\otimes_{i=1}^n V(Y)$ via the isomorphism (\ref{tp}) above and denote the
resulting linear operator by $\M_{i,j}^{t,l}$.
In Theorem \ref{tn} below we show that
the action of $\M_{i,j}^{t,l}$ on the basis 
$K_X(n)$ mirrors the action of $M_{i,j}^t$ on the standard basis of the
Boolean algebra $B(n)$.

Recall the maps $f_u$ on $V(X)$ from the introduction.
Define linear operatots $\Z, \R_u : V(Y)\rar V(Y),\;u=0,\ldots ,m,$ on $V(Y)$
as follows
\begin{itemize}
\item $\Z(L_0) = L_0$ and $\Z(x) = 0$ for $x\in X$.
\item For $u=0,\ldots ,m$, $\R_u(L_0) = 0$ and $\R_u(y) = f_u(y),\;y\in X$. 
\end{itemize}

Note that, for $w=0,\ldots ,m,$
\beq \label{diff1}
\R_u(v)&=&\lambda(u,w)v,\;v\in \B_w. 
\eeq

Let there be given a $(m+4)$-tuple 
$$\CS =(S_U,S_D,S_\Z,S_{\R_0},\ldots ,S_{\R_m})$$ 
of pairwise disjoint subsets of $[n]$ with
union $[n]$ (it is convenient to index the components of $\CS$ in this
fashion). 
Define a linear operator
$$ F(\CS): \bigotimes_{i=1}^n V(Y) \rar \bigotimes_{i=1}^n
V(Y) $$
by $F(\CS)=F_1 \otimes \cdots \otimes F_n$, where each $F_i$ is $U_1$ or
$D_1$ 
or $\Z$ or $\R_j$ according as $i \in S_U$ or $S_D$ or $S_\Z$
or $S_{\R_j}$, respectively.

Let $b\in \Bx(n)$. It follows from the definitions that
\beq \label{nz}
F(\CS)(\ob{b})\not= 0 &\mbox{ implies}& S_D \cup S_{\R_0} \cup \cdots 
\cup S_{\R_m} =
S(b),\;S_U \cup S_\Z = [n] - S(b).
\eeq

Given a $(m+4)$-tuple $r=(r_U,r_D,r_{\Z},r_{\R_0},\ldots ,r_{\R_m})\in
C(n,m+4)$ define
$\Pi(r)$ to be the set of all $(m+4)$-tuples 
$\CS =(S_U,S_D,S_\Z,S_{\R_0},\ldots ,S_{R_m})$ 
of pairwise
disjoint subsets of $[n]$ with union $[n]$ and with $|S_U|=r_U, |S_D|=r_D,
|S_\Z|= r_{\Z}$, and $|S_{\R_j}|=r_{\R_j}$, for $j=0,\ldots ,m$.

\bl \label{cbl} 
Let $(i,j,t,l)\in \I_X(n)$ with $l=(l_0,\ldots ,l_m)$.
Set $r=(i-t,j-t,n+t-i-j,l_0,\ldots ,l_m)$. Then
\beq \label{tt1}
\M_{i,j}^{t,l} &=& \sum_{\CS \in \Pi(r)} F(\CS).
\eeq 
\el
\pf 
Let $b= (b_1,\ldots ,b_n)\in \Bx(n)$ and $\CS=(S_U,S_D,S_\Z,S_{\R_0},\ldots
,S_{\R_m}) \in \Pi(r)$. We consider two cases:

(i) $|S(b)| \not= j$: In this case we have $\M_{i,j}^{t,l}
(\ob{b})=0$. Now $|S_D|+|S_{\R_0}|+\cdots +|S_{\R_m}|=j-t+l_0+\cdots +l_m=j$. 
Thus, from
(\ref{nz}), we also have $F(\CS)(\ob{b})=0$.  

(ii) $|S(b)|=j$: Assume that $F(\CS)(\ob{b})\not= 0$. We have from
(\ref{nz}) that $F(\CS)(\ob{b}) = \sum_{a} \ob{a}$,
where the sum is over all
$a=(a_1,\ldots ,a_n) \in \Bx(n)_i$
with $S(a) = S_U \cup S_{\R_0}\cup \cdots  \cup S_{\R_m}$
and $(a_k,b_k)\in Z_q$, for $k\in S_{\R_q},\;q=0,\ldots ,m$. 

Going over all elements of $\Pi(r)$ and summing we see that 
both sides
of (\ref{tt1}) evaluate to the same element on $\ob{b}$. $\Box$ 

\bt \label{tn}
Let $0\leq s \leq n$, $p=(p_1,\ldots ,p_m)\in C(s,m)$ and 
$(i,j,t,l)\in \I_X(n)$.
Let $(\A,f,B)\in \K_X(n,(n-s,p_1,\ldots ,p_m))$. Set $p_0=t-s$ and $p^+ =
(p_0,p_1,\ldots ,p_m)$.

(i) $\M_{i,j}^{t,l}(v(\A,f,B)) =0$ if $|B|\not= j-s$. 

(ii) If $|B| = j-s$ then 
\beq \label{cons}
\M_{i,j}^{t,l}(v(\A,f,B)) &=&  
       (|X|)^{\frac{i+j}{2}-t}\;\Lambda(\lambda,l,p^+)\;
 \left( {\ds{\sum_{B'}}}\;v(\A,f,B')\right), 
\eeq
where the sum is over all $B'\subseteq ([n]-\Sigma(\A))$ 
with $|B'|=i-s$ and $|B\cap B'| = t-s$.
\et
\pf
Let $r = (i-t,j-t,n+t-i-j,l_0,\ldots ,l_m)$, where $l=(l_0,\ldots ,l_m)$
 and let 
$\CS =
(S_U,S_D,S_\Z,S_{\R_0},\ldots ,S_{\R_m}) \in \Pi(r)$. 
Assume that $F(\CS)(v(\A,f,B))\not= 0$. 
Then we must have
\beq \label{ic1}
S_U \cup S_\Z &=& [n] - (\Sigma(\A) \cup B),\\ \label{ic2}
S_D\cup S_{\R_0} \cup \cdots \cup S_{\R_m} &=& \Sigma(\A) \cup B,\\
\label{ic3}
S_D &\subseteq& B.\eeq
Thus, using (\ref{ic2}) above, 
$|B|=j-t+l_0+\cdots +l_m - s=j-s$ (so part (i) follows).

We have
\beqn S_{\R_0} \cup \cdots \cup S_{\R_m} &=& \Sigma(\A) \cup (B-S_D).
\eeqn
Since each vector in $\B$ is an eigenvector for each of $\R_0,\ldots ,\R_m$
we see that
$$F(\CS)(v(\A,f,B))=\alpha \,v(\A,f,B'),$$
where $\alpha$ is a scalar and where $B'= S_U \cup (B-S_D)$, yielding$
|B'|=i-s$ and $|B\cap B'|= |B-S_D|=(j-s)-(j-t)=t-s$.

We now determine $\alpha$. Write $\A=(A_1,\ldots ,A_m)$ and put $A_0=B-S_D$.
It is easily seen, using (\ref{ipt3}) and the definition of $\lambda(u,w)$, that
\beq \label{ic4}
\alpha &=& |X|^{\frac{i+j}{2}-t}\;\prod_{u=0}^m \prod_{w=0}^m
\lambda(u,w)^{|S_{\R_u}\cap A_w|}.
\eeq 

The $\{0,1,\ldots ,m\}\times \{0,1,\ldots ,m\}$ 
integer matrix with entry in row $u$, column $w$
equal to $|S_{\R_u}\cap A_w|$ has row sums $l_0,\ldots ,l_m$ and columns
sums $p_0,p_1,\ldots ,p_m$.

Now fix a $\{0,1,\ldots ,m\}\times \{0,1,\ldots ,m\}$ 
integer matrix $(r(u,w))$ with 
row sums $l_0,\ldots ,l_m$ and columns sums $p_0,p_1,\ldots ,p_m$ and fix
$B'\subseteq ([n]-\Sigma(\A))$ with $|B'|=i-s$ and $|B\cap B'|=t-s$. 
We want to count the number of
$\CS =
(S_U,S_D,S_\Z,S_{\R_0},\ldots ,S_{\R_m}) \in \Pi(r)$
with $F(\CS)(v(\A,f,B)) = \alpha \,v(\A,f,B')$, where
$|S_{\R_u} \cap A_w|=r(u,w)$ for $u,w\in\{0,\ldots ,m\}$ and $\alpha$ is
given by (\ref{ic4}) above.

The discussion above shows that $S_U=B'-B$ and $S_D=B-B'$ and (\ref{ic1})
above then determines $S_{\Z}$. The number of choices for $S_{\R_0},\ldots
,S_{\R_m}$ is easily seen to be
$\left\{ \prod_{w=0}^m \binom{p_w}
{r(0,w),\ldots ,r(m,w)}\right\}$.

Going over all elements of $\Pi(r)$ and summing we get the result.
$\Box$

\noi {\bf Proof of Theorem \ref{ebdgba} (parts (iii) and (iv)):}

\noi (iii) Let $(i,j,t,l)\in \I_X(n)$ and 
$(k,s,p)\in \J_X(n)$. Observe that the term 
$$       (|X|)^{\frac{i+j}{2}-t}\;\Lambda(\lambda,l,p^+)$$
in (\ref{cons}) above depends only on $i,j,t,l,s,p$ (and not on $(\A, f, B)$).
It thus follows from Theorem \ref{tn} and Theorem
\ref{schrijver} (iii) that each SSJB in $J_X(n,k,s,p)$ is closed under all
the operators $\M_{i,j}^{t,l}$ and that all these SSJB's give rise to
identical blocks. That the image of $\Phi$ consists of all such block
diagonal matrices follows from 
the dimension count (\ref{card1}) and Lemma \ref{bi}.

\noi (iv) Follows from Theorem \ref{tn} and Theorem \ref{schrijver} (iv).
$\Box$

\bt \label{ev}
(Ceccherini-Silberstein,
Scarabotti, and Tolli {\bf\cite{cst1}})
Let $0\leq i \leq n$, $(t,l) \in \I_X(n,i)$, and $(k,s,p)\in \J_X(n,i)$ with
$p=(p_1,\ldots ,p_m)$. Set $p_0=t-s$ and $p^+=(p_0,p_1,\ldots ,p_m)$.
The eigenvalue of $M_{i,i}^{t,l}$ on $V_X(n,i,k,s,p)$ is
$$
(|X|)^{i -t}\;\Lambda(l,p^+)\;
\left\{ \sum_{u=0}^{n-s} 
(-1)^{u-t+s}\;\binom{u}{t-s}\binom{n-k-u}{i-s-u}\binom{i-k}{i-s-u}\right\}.
$$
\et
\pf Follows from
substituting $j=i$ in Theorem \ref{ebdgba} (iv) and noting that
$${\binom{n+s-2k}{i-k}}^{-1}{\binom{n+s-2k}{u+s-k}}{\binom{n-k-u}{i-s-u}} =
{\binom{i-k}{i-s-u}}.\;\;\;\;\Box$$

\end{document}